\documentclass[12pt]{article}
\begin{document}

\begin{center}
\bf Geometric angle structures on triangulated surfaces
\end{center}

\begin{center}
Ren Guo
\end{center}

\bigskip
\noindent {\bf Abstract} In this paper we characterize a function
defined on the set of edges of a triangulated surface such that there is a
spherical angle structure having the function as the edge invariant (or Delaunay invariant).
We also characterize a function such that there is a hyperbolic angle structure
having the function as the edge invariant.
\bigskip

\noindent \S 1. {\bf Introduction}
\medskip

Suppose $S$ is a closed surface and $T$ is a triangulation of $S$.
Here by a triangulation we mean the following: take a finite
collection of triangles and identify their edges in pairs by
homeomorphism. Let $V, E, F$ be the sets of all vertices, edges
and triangles in $T$ respectively. If $a, b$ are two simplices in
triangulation $T$, we use $a<b$ to denote that $a$ is a face of
$b$. Let $C(S,T)=\{ (e, f) | e \in E, f \in F,$ such that $e <
f\}$ be set of all \it corners \rm of the triangulation. An \it
angle structure \rm on a triangulated surface $(S,T)$ assigns each
corner of $(S,T)$ a number in $(0, \pi)$. A \it Euclidean (or
hyperbolic, or spherical) angles structure \rm is an angle
structure so that each triangle with the angle assignment is
Euclidean (or hyperbolic, or spherical). More precisely, a
Euclidean angle structure is a map $x: C(S,T)\to (0, \pi)$
assigning every corner $i$ (for simplicity of notation, we use one
letter to denote a corner) a positive number $x_i$ such that
$x_i+x_j+x_k=\pi$ whenever $i,j,k$ are three corners of a
triangle. A hyperbolic angle structure is a map $x: C(S,T)\to (0,
\pi)$ such that $x_i+x_j+x_k<\pi$. A spherical angle structure is
a map $x: C(S,T)\to (0, \pi)$ such that
\begin{equation} \label{1}
\left\{
\begin{array}{ccc}
x_i+x_j+x_k>\pi\\
x_j+x_k-x_i<\pi.
\end{array}
\right. \end{equation}

Actually it is proved in {\bf [B]} that positive numbers
$x_i,x_j,x_k$ are three inner angles of a spherical triangle if
and only if they satisfy conditions $(1)$.
\medskip

Given an angle structure $x: C(S,T)\to (0, \pi)$, we define its
\it edge invariant \rm which is a function $D_x: E \to (0,2\pi)$
such that $D_x(e)=x_i+x_{i'}$ where $i=(e,f),i'=(e,f')$ are two
opposite corners facing the edge $e$. And we define its \it
Delaunay invariant \rm which is a function $\mathcal{D}_x: E \to
(-2\pi,2\pi)$ such that
$\mathcal{D}_x(e)=x_j+x_k+x_{j'}+x_{k'}-x_i-x_{i'}$ where
$i=(e,f),i'=(e,f')$ are two opposite corners facing the edge $e$
and $j,k$(or $j',k'$) are the other two corners of the triangle
$f$ (or $f'$).
\medskip

For the simplicity of natation, we use $G$ to denote a fixed geometry,
where $G=E,H$ or $S$ means the Euclidean, hyperbolic or spherical geometry
respectively. Now given a function $D: E \to (0,2\pi)$ (or
$\mathcal{D}: E \to (-2\pi,2\pi)$), we use $AG(S,T;D)$ (or $AG(S,T;\mathcal{D})$)
to denote the set fo all $G$ angle structures having $D$ (or
$\mathcal{D}$) as the edge (or Delaunay) invariant.
\medskip

The motivation of considering these sets is the study of \it
geometric cone metrics \rm with prescribed edge invariant or
Delaunay invariant on triangulated surfaces from the variational point of view.
A \it Euclidean (or
hyperbolic, or spherical) cone metric \rm assigns each edge in $T$
a positive number such that the numbers on any three edges of a triangle in $T$
form three edge length of a Euclidean (or
hyperbolic, or spherical) triangle. The variational method contains a variational problem
and a linear programming problem. The variational problem is to
show that the unique maximal point of a convex "capacity" defined on the
set $AG(S,T;D)$ (or $AG(S,T;\mathcal{D})$) gives the unique geometric
cone metric. The linear programming problem is to characterize the function $D$ (or
$\mathcal{D}$) such that the set $AG(S,T;D)$ (or $AG(S,T;\mathcal{D})$) is nonempty.
\medskip

For Euclidean angle
structures, the Delaunay invariant and the edge invariant are
related by $2D_x(e)+\mathcal{D}_x(e)=2\pi$ for any $e$. Thus given
two functions $D$ and $\mathcal{D}$ satisfying
$2D(e)+\mathcal{D}(e)=2\pi$ for any $e$, we have
$AE(S,T;D)=AE(S,T;\mathcal{D}).$ Therefore the problem of
Euclidean cone metric with given edge invariant is eqivalent to
the problem of Euclidean cone metric with given Delaunay invariant.
Rivin {\bf[Ri1]} {\bf[Ri2]} worked out the variational problem and the
linear programming problem about $AE(S,T;D).$
Leibon {\bf[Le]} worked out the variational problem and the
linear programming problem about $AH(S,T;\mathcal{D}).$
Luo {\bf[Lu]} worked out the variational problem about $AS(S,T;D)$
the linear programming problem about which will be solved in this paper (theorem 1).
Although the variational problems about $AH(S,T;D)$ and $AS(S,T;\mathcal{D})$ are still open,
we will solve the linear programming problem about them in this paper (theorem 2 and 3).
\medskip

The main results are the following.
For a triangulated surface $(S,T),$ a subset $X\subseteq F,$ we
use $|X|$ to denote the number of triangles in $X$ and we use
$E(X)$ to denote the set of all edges of triangles in $X.$
\medskip

\noindent {\bf Theorem 1.} \it Given a triangulated surface
$(S,T)$ and a function $D: E\to (0,\pi)$, the set $AS(S,T;D)$ is
nonempty if and only if for any subset $X\subseteq F,$
$$\pi |X|< \sum _{e\in E(X)}D(e).$$ \rm
\medskip

\noindent {\bf Theorem 2.} \it Given a triangulated surface
$(S,T)$ and a function $D: E\to (0,2\pi)$, the set $AH(S,T;D)$ is
nonempty if and only if for any subset $X\subset F,$
$$\pi(|F|-|X|)> \sum _{e\notin E(X)}D(e).$$ \rm
\medskip

\noindent {\bf Theorem 3.} \it Given a triangulated surface
$(S,T)$ and a function $\mathcal{D}: E\to (-2\pi,2\pi)$, the set
$AS(S,T;\mathcal{D})$ is nonempty if and only if for any subset
$X\subset F,$
$$\pi(|F|-|X|)> \sum _{e\notin E(X)}(\pi-\frac12\mathcal{D}(e)).$$ \rm
\medskip

The paper is organized as follows. In section 2, we prove theorem
1 by using Leibon's result. In section 3, we recall the duality
theorem in linear programming. In section 4, following Rivin's
method, we prove theorem 2 and 3 by using the duality theorem.
\medskip

\noindent{\bf Acknowledgement} I wish to thank my advisor,
Professor Feng Luo, for suggesting this problem and for fruitful
discussion.
\medskip

\noindent \S 2. {\bf Prove of theorem 1}
\medskip

First let us recall the Leibon's result of characterization of the
function $\mathcal{D}$ such that the set $AH(S,T;\mathcal{D})$ is
nonempty.
\medskip

\noindent {\bf Theorem 4.}(Leibon){\bf[Le]} \it Given a
triangulated surface $(S,T)$ and a function $\mathcal{D}: E\to
(0,2\pi)$, the set $AH(S,T;\mathcal{D})$ is nonempty if and only if
for any subset $X\subseteq F,$
$$\pi |X|< \sum _{e\in E(X)}(\pi-\frac12\mathcal{D}(e)).$$ \rm
\medskip

\noindent {\bf Proof of theorem 1.} To show the conditions are
necessary, for any  $X\subseteq F$, we have $\sum _{e\in
E(X)}D(e)=\sum_{e\in E(X)}(x_i+x_{i'}),$ where $i,i'$ are two
opposite corners facing the edge $e$. It turns out that the right
hand side of the equation is equal to $\sum_{f\in
X}(x_i+x_j+x_k)+\sum x_h,$ where the corner $h=(e,f^*)$ with $e\in
E(X)$ and $f^*\notin X.$ Hence $\sum _{e\in E(X)}D(e)\geq
\sum_{f\in X}(x_i+x_j+x_k)> \sum_{f\in X}\pi= \pi|X|.$
\medskip

To show the conditions are sufficient, let us define a function
$\mathcal{D}: E \to (0,2\pi)$ by setting
$\mathcal{D}(e)=2\pi-2D(e).$ Thus the conditions $\pi |X|< \sum
_{e\in E(X)}D(e)$ are equivalent to $\pi |X|< \sum _{e\in
E(X)}(\pi-\frac12\mathcal{D}(e))$ which guarantee
$AH(S,T;\mathcal{D})$ is nonempty by theorem 4. It follows that
there is a solution for the inequalities
$$
\left\{
\begin{array}{ccc}
x_i+x_j+x_k< \pi & \ i,j,k\ \mbox{are three corners of a triangle}\\
x_j+x_k+x_{j'}+x_{k'}-x_i-x_{i'}=\mathcal{D}(e) \\
x_i> 0
\end{array}
\right.
$$
\medskip

Let us define new variables $y_i$ for all $i \in C(S,T)$ by setting
$$y_i=\frac{\pi+x_i-x_j-x_k}{2}$$ provided $i,j,k$ are three corners
of a triangle. And since $\mathcal{D}(e)=2\pi-2D(e),$
the inequalities above are equivalent to
$$
\left\{
\begin{array}{ccc}
y_i+y_j+y_k> \pi& \ i,j,k\ \mbox{are three corners of a triangle}\\
y_i+y_{i'}=D(e)& \ i,i'\ \mbox{are two
opposite corners facing an edge}\ e \\
y_j+y_k< \pi &  \ j,k\ \mbox{are two corners of
a triangle}\\
\end{array}
\right.
$$

This solution obviously satisfies
$$\left\{
\begin{array}{ccc}
y_i+y_j+y_k> \pi& \ i,j,k\ \mbox{are three corners of a triangle}\\
y_i+y_{i'}=D(e)& \ i,i'\ \mbox{are two
opposite corners facing an edge}\ e \\
y_j+y_k-y_i< \pi &  \ i,j,k\ \mbox{are three corners
of a triangle}\\
y_i> 0\\
\end{array}
\right.
$$
Thus we obtain an angle structure in $AS(S,T;D)$. QED
\medskip

\noindent \S 3. {\bf Duality Theorem}
\medskip

We fix the notations as follows: $x=(x_1, ..., x_n)^t$ is a column
vector in $ \mathbf{R}^n$. The standard inner product in $
\mathbf{R}^n$ is denoted by $a^t x$. If  $A: \mathbf{R}^n \to
\mathbf{R}^m$ is a linear transformation, we denote its transpose
by $A^t: \mathbf{R}^m \to \mathbf{R}^n.$ Given two vectors $x, a$
in $ \mathbf{R}^n$, we say $x \geq a$ if $x_i \geq a_i$ for all
indices $i$. Also $x>a$ means $x_i > a_i$ for all indices $i$.
\medskip

A linear programming problem $(P)$ is to minimize an \it objective
function \rm $z = a^tx$ subject to the \it restrain conditions\rm
$$
\left\{
\begin{array}{ccc}
Ax=b\\
x \geq 0
\end{array}
\right.
$$ where $x \in \mathbf{R}^n$, $b\in \mathbf{R}^m$ and $A:
\mathbf{R}^n \to \mathbf{R}^m$ is a linear transformation. We call
a point $x$ satisfying the restrain conditions a \it feasible
solution \rm and denote the set of all the feasible solutions by
$D(P)=\{x \in \mathbf{R}^n|Ax=b,x \geq 0\}.$ An \it optimal
solution \rm $x$ for $(P)$ is a feasible solution so that the
objective function $z$ realizes the minimal value. The \it dual
problem \rm $(P^*)$ of $(P)$ is to maximize $z = b^t y$ subject to
$ A^t y \leq a , y\in \mathbf{R}^m$. Let us recall the duality
theorem in linear programming. The proof of the theorem can be
found in the book {\bf[KB]}.

\medskip
\noindent {\bf Theorem 5.}  \it The following statements are
equivalent.

\noindent (a) Problem (P) has an optimal solution.

\noindent (b) $D(P) \neq \emptyset$ and $D(P^*) \neq \emptyset$.

\noindent (c) Both problem $(P)$ and problem $(P^*)$ have optimal solutions so
that the minimal value of $(P)$ is equal to the maximal value of
$(P^*)$. \rm
\medskip

In applications that we are interested, there is a special case
that the objective function $z =0$ for $(P)$. Thus the optimal
solution exists if and only if $D(P) \neq \emptyset$. Thus we
obtain the following corollary.
\medskip

\noindent {\bf  Corollary 6.} For $A: \mathbf{R}^n \to
\mathbf{R}^m$ and $b\in \mathbf{R}^m,$ the set $\{ x \in
\mathbf{R}^n |Ax=b, x \geq 0\} \neq \emptyset$  if and only if the
maximal value of $z = b^ty$ on  $\{ y \in R^m | A^t y \leq 0\}$ is
non-positive.
\medskip

\noindent \S 4. {\bf Proof of theorem 2 and 3 }
\medskip

By following Rivin's method in {\bf [Ri2]}, we will prove a lemma
about the closure of $AH(S,T;D)$ in $ \mathbf{R}^{3|F|}=\{(x_i)^t, i\in C(S,T)\}$. The
closure of $AH(S,T;D)$ consists of all the points satisfying
$$
\left\{
\begin{array}{ccc}
x_i+x_j+x_k \leq \pi& \ i,j,k\ \mbox{are three corners of a triangle}\\
x_i+x_{i'}=D(e) & \ i,i'\ \mbox{are two opposite corners facing an edge}\ e\\
x_i\geq 0 \\
\end{array}
\right.
$$
\medskip

\noindent {\bf Lemma 7.} Given a triangulated surface $(S,T)$ and
a function $D: E\to [0,2\pi]$, the closure of $AH(S,T;D)$ is
nonempty if and only if for any subset $X\subset F,$
$$\pi(|F|-|X|)\geq \sum _{e\notin E(X)}D(e).$$ \rm
\medskip

\noindent {\bf Proof.} The linear programming problem $(P)$ with
variables $x=(...,x_i,...,t_f,...)$ indexed by $C(S,T)\cup F$ is
to minimize the objective function $z = 0$ subject to the restrain
conditions
$$
\left\{
\begin{array}{ccc}
x_i+x_j+x_k+ t_f=\pi& \ i,j,k\ \mbox{are three corners of a triangle} f\\
x_i+x_{i'}=D(e)& \ i,i'\ \mbox{are two opposite
corners facing an edge}\ e\\
x_i\geq 0 \\
t_f\geq 0
\end{array}
\right.
$$
The dual problem $(P^*)$ with variable $y =( ...,y_f, ..., y_e,
...)$ indexed by $E \cup F$ is to maximize the objective function
$z = \sum_{f \in F} \pi y_f+ \sum_{e \in E} D(e) y_e$ subject to
the restrain conditions
$$
\left \{
\begin{array}{ccc}
y_f \leq 0&\\
y_f + y_e \leq 0&\  \mbox{whenever}\ e < f.
\end{array}
\right.
$$
Since the closure of $AH(S,T;D)$ is nonempty is equivalent to that the set
$D(P)$ is nonempty, by corollary 6, the latter one is equivalent
to that the maximal value of the objective function of $(P^*)$ is
non-positive.
\medskip

To show the conditions $\pi(|F|-|X|)\geq \sum _{e\notin E(X)}D(e)$
for any $X\subset F$ are necessary, for any $X\subset F,$ let
$$y_f=
\left\{
\begin{array}{ccc}
0&\mbox{if}\ f\in X \\
-1&\mbox{if}\ f\notin X
\end{array}
\right.\     \mbox{and}\     y_e= \left\{
\begin{array}{ccc}
0&\mbox{if}\ e\in E(X) \\
1&\mbox{if}\ e\notin E(X)
\end{array}
\right.
$$

We claim that $(y_f, y_e)$ is a feasible solution. In fact, given
a pair $e<f,$ if $f\in X$, we must have $e\in E(X)$, then $y_f +
y_e=0.$ If $f\notin X$, then $y_f + y_e=-1+y_e \leq 0$.
\medskip

By the assumption that the maximal value of the objective function
of $(P^*)$ is non-positive, since $(y_f, y_e)$ is feasible, we
have $0\geq z(y_f, y_e) = \sum_{f \notin X} \pi y_f+ \sum_{e
\notin E(X)} D(e) y_e = \pi(|X|-|F|)+ \sum _{e\notin E(X)}D(e).$
\medskip

To show the conditions are sufficient, take an arbitrary feasible
solution $(y_f, y_e)$. If $y_f=0$ for all $f$, from $y_f+y_e\leq
0$, we know $y_e\leq 0$. Hence $z(y_f,y_e) = \sum _{e\notin
E}D(e)y_e \leq 0$, since $D(e)\in [0,2\pi].$ Otherwise, define
$X=\{f\in F | y_f = 0\}\subset F$, and let $a=max\{y_f, f\notin
X\}$. We have $a < 0$. Define
$$y_f^{(1)}=
\left\{
\begin{array}{ccc}
y_f=0&\mbox{if}\ f\in X \\
y_f-a&\mbox{if}\ f\notin X
\end{array}
\right.\     \mbox{and}\     y_e^{(1)}= \left\{
\begin{array}{ccc}
y_e&\mbox{if}\ e\in E(X) \\
y_e+a&\mbox{if}\ e\notin E(X)
\end{array}
\right.
$$

We claim that $(y_f^{(1)},y_e^{(1)})$ is a feasible solution. In
fact, $y_f^{(1)}\leq 0$. Given a pair $e<f,$ if $f\in X$, we must
have $e\in E(X)$, then $y_f^{(1)}+y_e^{(1)}=y_f+y_e\leq 0$. If $f\notin X$
and $e\notin E(X)$, then $y_f^{(1)}+y_e^{(1)}=y_f-a+y_e+a\leq 0$. If $f\notin
X$ but $e\in E(X)$, there exists another triangle $f'\in X$ so
that $e<f'$, then $y_e=y_e+y_{f'}\leq 0$. Therefore
$y_f^{(1)}+y_e^{(1)} = y_f-a+y_e\leq y_f-a \leq 0$, since $a$ is
the maximum.
\medskip

Now the value of the objective function is
$z(y_f^{(1)},y_e^{(1)})= z(y_f,y_e)+ a(\pi (|X|-|F|)+\sum
_{e\notin E(X)}D(e))\geq z(y_f,y_e)$, according to the conditions.
Note the number of 0's in $\{y_f^{(1)}\}$ is more than that in
$\{y_f\}$. By the same procedure, after finite steps, it ends at a
feasible solution $(y_f^{(n)}=0,y_e^{(n)})$. We have
$z(y_f^{(n)},y_e^{(n)})\leq 0$. Since the value of the objective
function does not increase, therefore $0\geq
z(y_f^{(n)},y_e^{(n)})\geq \ldots \geq z(y_f^{(1)},y_e^{(1)}) \geq
z(y_f, y_e)$. QED
\medskip

\noindent {\bf  Proof of theorem 2.} Let $x_i=a_i+\varepsilon$ for any $i \in C(S,T),$
where $a_i\geq 0$ and $\varepsilon\geq 0.$ The linear programming
problem $(P)$ with variables $\{...,a_i,...\varepsilon\}$ is to
minimize the objective function $z = -\varepsilon$ subject to the
restrain conditions
$$
\left\{
\begin{array}{ccc}
a_i+a_j+a_k+3\varepsilon\leq\pi& \ i,j,k\ \mbox{are three corners of a triangle}\\
a_i+a_j+2\varepsilon=D(e)& \ i,j\ \mbox{are two
opposite corners facing an edge}\ e \\
a_i\geq 0\\
\varepsilon \geq 0
\end{array}
\right.
$$
The dual problem $(P^*)$ with variable $y =( ...,y_f, ..., y_e,
...)$ indexed by $E \cup F$ is to maximize the objective function
$z = \sum_{f \in F} \pi y_f+ \sum_{e \in E} D(e) y_e$ subject to
the restrain conditions
$$
\left\{
\begin{array}{ccc}
y_f \leq 0\\
y_f + y_e \leq 0&  \mbox{whenever} e < f\\
3\sum_{f \in F} y_f + 2 \sum_{e \in E} y_e \leq -1
\end{array}
\right.
$$
By the theorem 5(c), the maximal value of the objective function
of $(P^*)$ is negative is equivalent to that the minimal value of
the objective function of $(P)$ is negative. The latter one is
equivalent to that there exists a feasible solution $a_i\geq 0,
\varepsilon> 0$. Therefore the set $AH(S,T;D)$ is nonempty.
\medskip

We only need to show that the maximal value of the objective
function of $(P^*)$ is negative is equivalent to the conditions
$\pi(|F|-|X|)> \sum _{e\notin E(X)}D(e)$ for any $X\subset F.$
\medskip

To show the conditions are necessary, for any $X\subset F,$ we
have $2|E(X)|>3|X|$ or $2|E(X)|\geq 3|X|+1$. Let
$$y_f=
\left\{
\begin{array}{ccc}
0&\mbox{if}\ f\in X \\
-1&\mbox{if}\ f\notin X
\end{array}
\right.\     \mbox{and}\     y_e= \left\{
\begin{array}{ccc}
0&\mbox{if}\ e\in E(X) \\
1&\mbox{if}\ e\notin E(X)
\end{array}
\right.
$$
We claim that $(y_f, y_e)$ is a feasible solution. If fact, as in
lemma 7, we can check $y_f+ y_e\leq 0$ for any pair $e<f.$
Furthermore $$3\sum_{f \in F} y_f + 2 \sum_{e \in E} y_e
=3\sum_{f\notin X}(-1)+2\sum_{e\notin
E(X)}1=3(|X|-|F|)+2(|E|-|E(X)|)$$
$$=3|X|-2|E(X)|+2|E|-3|F|=3|X|-2|E(X)|\leq-1$$
since $2|E|=3|F|$. Now $(y_f,y_e)$ is feasible
implies that $z(y_f,y_e)<0$ which is equivalent to $\pi(|F|-|X|)<
\sum _{e\notin E(X)}D(e).$
\medskip

To show the conditions are sufficient, by the proof of lemma 7 we
know the maximal value of the objective function of $(P^*)$ is
$\leq 0$ under the conditions. We try to show it can not be 0.
Assume that $(y_f,y_e)$ is a feasible solution satisfying
$z(y_f,y_e)=0.$ We claim that $y_f=0$ for all $f$. Otherwise, as
in the proof of lemma 7, we can find another feasible solution
$(y_f^{(1)},y_e^{(1)})$ and we can check that
$z(y_f^{(1)},y_e^{(1)}) = z(y_f,y_e)+ a(\pi (|X|-|F|)+\sum
_{e\notin E(X)}D(e)) > z(y_f,y_e)= 0$, according to the
conditions. It is contradiction since the maximal value of the
objective function of $(P^*)$ is $\leq 0.$
\medskip

Now from $y_f=0$ for all $f$ we see $y_e\leq 0$. Since $0=z(y_f,y_e)=\sum_{e
\in E} D(e) y_e$ and $D(e)>0$, we get $y_e=0$ for all $e$ and therefore $(y_f,y_e)= (0,0).$
But $(y_f,y_e)= (0,0)$ does not satisfy $3\sum_{f \in F} y_f + 2
\sum_{e \in E} y_e \leq -1.$ It is a contradiction since we assume
that $(y_f,y_e)$ is a feasible solution. This proves that the
maximal value of the objective function of $(P^*)$ is negative.
QED
\medskip

\noindent {\bf  Proof of theorem 3.} Given two functions
$D:E\to (0, 2\pi)$ and $\mathcal{D}: E\to (-2\pi,2\pi)$ satisfying
$2D(e)+\mathcal{D}(e)=2\pi$ for any $e$, we claim that
$AH(S,T;D)\neq \emptyset$ is eqivalent to $AS(S,T;\mathcal{D})\neq \emptyset.$
By this claim, theorem 3 is true
as a corollary of theorem 2.
\medskip

In fact, $AS(S,T;\mathcal{D})$ is the set of solutions for the inequalities
$$
\left\{
\begin{array}{ccc}
x_i+x_j+x_k > \pi& \ i,j,k\ \mbox{are three corners of a triangle}\\
x_j+x_k-x_i < \pi & \ i,j,k\ \mbox{are three corners of a triangle}\\
x_j+x_k+x_{j'}+x_{k'}-x_i-x_{i'}=\mathcal{D}(e) \\
x_i > 0
\end{array}
\right.
$$

Let us define new variables $y_i$ for all $i\in C(S,T)$ by setting
$$y_i=\frac{\pi+x_i-x_j-x_k}{2}$$ provided $i,j,k$ are three corners
of a triangle. Since $2D(e)+\mathcal{D}(e)=2\pi,$ we see that the
inequalities above are equivalent to
$$
\left\{
\begin{array}{ccc}
y_i+y_j+y_k < \pi& \ i,j,k\ \mbox{are three corners of a triangle}\\
y_i> 0 \\
y_i+y_{i'}=D(e) & \ i,i'\ \mbox{are two opposite corners facing an edge}\ e\\
y_j+y_k<\pi & \ j,k\ \mbox{are two corners of a triangle}\\
\end{array}
\right.
$$

Since $y_i+y_j+y_k < \pi$ implies $y_j+y_k<\pi$, we can omit the
latter one. Equivalently, we get
$$
\left\{
\begin{array}{ccc}
y_i+y_j+y_k < \pi& \ i,j,k\ \mbox{are three corners of a triangle}\\
y_i> 0 \\
y_i+y_{i'}=D(e) & \ i,i'\ \mbox{are two opposite corners facing an edge}\ e\\
\end{array}
\right.
$$

Now the set of solutions of the inequalities above is exactly $AH(S,T;D).$
Thus we see $AH(S,T;D)\neq \emptyset$ is eqivalent to $AS(S,T;\mathcal{D})\neq \emptyset.$  QED
\medskip

\noindent {\bf Reference}
\medskip

[B] Marcel Berger, Geometry II. Springer-Verlag 1987

[KB]Bernard Kolman $\&$ Robert Beck,  Elementary
Linear Programming with Applications. Academic Press 2 edition
1995

[Le] Gregory Leibon,  Characterizing the Delaunay
decompositions of compact hyperbolic surface. Geom. Topol.
6(2002), 361-391

[Lu] Feng Luo,  A Characterization of spherical
polyhedron surfaces.\\
http://front.math.ucdavis.edu/math.GT/0408112

[Ri1] Igor Rivin,  Euclidean structures on simplicial
surfaces and hyperbolic volume. Ann. of Math. (2) 139 (1994), no.
3, 553-580

[Ri2] Igor Rivin,  Combinational optimization in
geometry. Advance in Applied Math. 31(2003), no. 1, 242-271

\bigskip
\noindent Department of Mathematics\\ Rutgers University\\
Piscataway, NJ 08854, USA

\bigskip
\noindent Email: renguo$@$math.rutgers.edu
\end{document}